\documentclass[12pt]{article}
\usepackage{amsxtra,amssymb,amsthm,amsmath,latexsym}

\theoremstyle{plain}

\newtheorem{theorem}{Theorem}[section]

\newtheorem{remark}{Remark}[section]
\newtheorem*{assumption}{Assumption}
\numberwithin{equation}{section}

\newcommand{\refT}[1]{Theorem~\ref{T:#1}}

\def\oH{\buildrel\circ\over H}
\def\oH1{\buildrel\circ\over H\kern-.02in{}^1}

\def\const{\hbox{\,const\,}}
\def\be{\begin{equation}}
\def\ee{\end{equation}}

\begin{document}
\title{                  
Dynamical systems method (DSM) for unbounded operators}

\author{
A.G. Ramm\\
 Mathematics Department, Kansas State University, \\
 Manhattan, KS 66506-2602, USA\\
ramm@math.ksu.edu\\
http://www.math.ksu.edu/\,$\widetilde{\ }$\,ramm}

\date{}

\maketitle\thispagestyle{empty}

\begin{abstract}
\footnote{Math subject classification:
35R25, 35R30, 37B55, 47H20, 47J05, 49N45, 65M32, 65R30}
\footnote{key words:
dynamical systems method, nonlinear operator equations,
ill-posed problems}

Let $L$ be an unbounded linear operator
in a real Hilbert space $H$, a generator
of $C_0$ semigroup, and $g:H\to H$ be a $C^2_{loc}$
nonlinear map. The DSM (dynamical systems method) for solving equation
$F(v):=Lv+gv=0$ consists of solving the Cauchy problem $\dot 
{u}=\Phi(t,u)$,
$u(0)=u_0$, where $\Phi$ is a suitable operator, and proving that
i) $\exists u(t) \quad \forall t>0$,
ii) $\exists u(\infty)$,
and iii) $F(u(\infty))=0$.

Conditions on $L$ and $g$ are given which allow one to choose $\Phi$ 
such that i), ii), and iii) hold.

\end{abstract}  

\section{Introduction}
Let $H$ be a real Hilbert space, $L$ be a linear, densely defined in $H$,
closed operator, a generator of $C_0$ semigroup (see\cite{P}),
$g:H\to H$ be a nonlinear $C^2_{loc}$ map, i.e.,
\be \label{e1.1}
 \sup_{u\in B(u_0,R)}  \| g^{(j)} (u)\| \leq m_j,
 \quad j=0,1,2, \quad B(u_0,R):=\{u:\|u-u_0\|\leq R\}, \ee
where $g^{(j)}$ is the Fr\'echet derivative of order $j$,
$R>0$ is some number, and $u_0\in H$ is some element.
In many applications the problems can be formulated as the following
operator equation:
\be \label{e1.2} F(v):=Lv+g(v)=0. \ee
We want to study this equation by the dynamical sysstems method (DSM),
which allows one also to develop numerical methods for solving 
equation (1.2).
The DSM for solving equation (1.2) consists of solving the problem:
\be \label{e1.3} \dot u= \Phi (t,u),\quad u(0)=u_0, \ee
where $\dot {u}:=\frac{du}{dt}$, and
 $\Phi(t,u)$ is a nonlinear operator chosen so that problem (1.3)
has a unique global solution which stabilizes at infinity to the solution
of equation (1.2):
\be \label{e1.4}
 \hbox{i)\ } \exists u(t) \forall t>0,
 \qquad \hbox{ii)\ } \exists u(\infty),
 \qquad \hbox{iii)\ } F(u(\infty))=0. \ee

In [2] the DSM has been studied and justified for
$F\in C^2_{loc}$ and
\be \label{e1.5}
 \sup_{u\in B(u_0,R)} \| [F'(u)]^{-1}\| \leq m_1; \ee
for monotone $F\in C^2_{loc}$;
for monotone hemicontinuous defined on all of $H$ operators $F$; 
for non-monotone $F\in C^2_{loc}$
such that there exists a $y$ such that $F(y)=0$ and the operator 
$A:=F'(y)$ maps any ball $B(0,r)$ centered at the origin and of 
sufficiently small
radius $r>0$ into a set which has a non-empty intersection with $B(0,R)$; 
and for $F\in C^2_{loc}$
satisfying a spectral condition:
$\|(F'(u)+\varepsilon)^{-1}\| \leq (c\varepsilon)^{-1}$,
$0< c\leq 1$, in which case $F$ is replaced by $F+\varepsilon I$
in \eqref{e1.3}, and then $\varepsilon$ is taken to zero.

{\it In this paper the DSM is justified for a class of nonlinear
unbounded operators of the type $L+g$, where $L$ is a generator of a $C_0$
semigroup, $g\in C^2_{loc}$, and some suitable additional assumptions
are made.}

Which assumptions are suitable?
A simple example is:
\be \label{e1.6}  \|L^{-1}\|\leq m. \ee
If \eqref{e1.6} holds, then \eqref{e1.2} is equivalent to
\be \label{e1.7} f(u):=u+L^{-1} g(u)=0, \quad f\in C^2_{loc}. \ee
Assume
\be \label{e1.8}
 \sup_{u\in B(u_0,R)} \|[I+L^{-1} g'(u)]^{-1} \| \leq m_1. \ee
This assumption, holds, e.g., if $L^{-1} g'(u)$ is a compact 
operator
in $H$ for any $u\in B(u_0,R)$ and the operator $I+L^{-1} g'(u)$
is injective.

Our first result is the following theorem:

\begin{theorem} \label{T:1.1}
Assume \eqref{e1.6}, \eqref{e1.8}, and let
$\Phi:=-[I+L^{-1}g'(u)]^{-1} [u+L^{-1} g(u)]$.
If
\be \label{e1.9} \|u_0+L^{-1} g(u_0)\| m_1\leq R, \ee
then equation \eqref{e1.2} has a unique solution
$v\in B(u_0,R)$, the conclusions
i), ii), and iii) hold for \eqref{e1.3}, and $u(\infty)=v$.
\end{theorem}

\begin{remark}\label{R:1.2}
If $L$ is not boundedly invertible, i.e., \eqref{e1.6} fails,
then one can use the following assumption (A):

\begin{assumption}[A]
{\it There exists a sector $S=\{z:0<|z|\leq a$,
$|\arg z-\pi|\leq \delta\}$, which consists of regular points of $L$.
Here $a>0$ and $\delta>0$ are arbitrary small fixed numbers.}
\end{assumption}
\end{remark}
If assumption (A) holds, then
\be \label{e1.10}
 \|(L+\varepsilon)^{-1}\| \leq \frac{1}{\varepsilon \sin (\delta)}. \ee
This estimate holds in particular if $L=L^*\geq 0$, and in this case
$ \sin (\delta)=1$.

 The following theorem is our next result: 

\begin{theorem}\label{T:1.2}
Assume that $L^\ast=L\geq 0$ is a densely defined linear operator,
$g\in C^2_{loc}$, $g'(u)\geq 0$ $\forall u\in H$,
equation \eqref{e1.2} is solvable and $v$ is its (unique)
minimal-norm solution.  Define
$\Phi(u)=-[I+(L+\varepsilon)^{-1}g'(u)]^{-1}
[u+(L+\varepsilon)^{-1}g(u)]$, $\varepsilon=\const>0$.
Assume that (1.8) holds with $L_\varepsilon:=L+\varepsilon I$
replacing $L$, and $m_1=m_1(\varepsilon)>0$. 
Then problem \eqref{e1.3} has a unique global solution
$u_\varepsilon(t)$, there exists $u_\varepsilon(\infty):=v_\varepsilon$,
and $F(v_\varepsilon):=L v_\varepsilon + g(v_\varepsilon)=0$.
Moreover, there exists the limit $\lim_{\varepsilon\to 
0}v_\varepsilon=v$,
which is the unique minimal-norm solution to \eqref{e1.2}.
\end{theorem}

In section 2 proofs are given.

\section{Proofs}

\begin{proof}[Proof of Theorem 1.1]
If $\Phi=-[I+L^{-1} g'(u)]^{-1} [u+L^{-1}g(u)]$,
and $p(t):=\|u+L^{-1}g(u)\|$, then $p\dot p=-p^2$.
Thus
\be \label{e2.1} p(t)=p(0) e^{-t}. \ee
From \eqref{e2.1}, \eqref{e1.8} and \eqref{e1.3} one gets
\be \label{e2.2}
\|\dot u\| \leq m_1 p(0)e^{-t}, \quad p(0)=\|u_0+L^{-1}g(u_0)\|. \ee
Inequality \eqref{e2.2} implies the global existence of 
$u(t)$,  the existence of
$u(\infty):=\lim_{t\to\infty} u(t)$, and the estimates:
\be \label{e2.3}
 \|u(t)-u(\infty)\| \leq m_1 p(0)e^{-t},
 \quad \|u(t)-u_0\| \leq m_1 p(0). \ee
If \eqref{e1.9} holds, then \eqref{e2.3} implies
$\|u(t)-u_0\|\leq R$, so the trajectory $u(t)$ stays in the ball
$B(u_0,R)$, that is, $ u(t)\in B(u_0,R)$ $\forall t\geq 0$.
Passing to the limit $t\to\infty$ in equation \eqref{e1.3} yields
\be \label{e2.4}
 0=-[I+L^{-1} g'(u(\infty))]^{-1}
 [u(\infty)+L^{-1} g(u(\infty))]. \ee
Thus $u(\infty):=v$ solves the equation $v+L^{-1}g(v)=0$, so $v$ solves 
(1.2), and therefore
i), ii) and iii) hold. \refT{1.1} is proved.
\end{proof}

{\it Proof of Theorem 1.2.}
If  $L^*=L\geq 0$ in $H$ and
$g'(u)\geq 0$, 
then, for any $\varepsilon>0$, \refT{1.1} yields the existence
of a unique solution $v_\varepsilon$ to equation \eqref{e1.2} with $L$
replaced by $L+\varepsilon I$. This solution 
$v_\varepsilon=u_\varepsilon(\infty)$,
where $u_\varepsilon(t)$ is the solution to \eqref{e1.3}
with 
$$\Phi=-[I+(L+\varepsilon)^{-1}g'(u)]^{-1}
[u+(L+\varepsilon)^{-1}g(u)].$$
Let us prove that $\lim_{\varepsilon\to 0} v_\varepsilon=v$,
where $v$ solves \eqref{e1.2}. We do not assume that 
\eqref{e1.2} has a unique solution.

Let $v_\varepsilon-v:=w$.
Then $Lw+\varepsilon v_\varepsilon + g(v_\varepsilon)-g(v)=0$,
so, using the assumptions $L\geq 0$ and $g'(u)\geq 0$, one gets
$\varepsilon(v_\varepsilon, v_\varepsilon-v)\leq 0$,
$\|v_\varepsilon\|^2\leq \|v_\varepsilon\| \|v\|$,
and $\|v_\varepsilon\|\leq \|v\|$, $\forall \varepsilon>0$.
Thus $v_\varepsilon\rightharpoonup v_0$ and 
$Lv_\varepsilon+g(v_\varepsilon)\to 0$,
where $\rightharpoonup$ stands for the weak convergence in $H$ and the
convergent subsequence is denoted $v_\varepsilon$ again.

In the above argument the element $v$ can be an arbitrary
element in the set $ N_F:=\{v: Lv+g(v)=0\}$. Thus, we have proved 
that $||v_0||\leq ||v||$ for all $v\in N_F$.

{\it Let us prove that $L(v_0)+g(v_0)=0$, i.e., $v_0\in N_F$.} 

Assume 
first that $v_0\in D(L)$. We
prove this assumption later. The monotonicity of $L+g$ yields:
\be \label{e2.5}
 (L(v_\varepsilon)+g(v_\varepsilon)+\varepsilon v_\varepsilon
 -L(v_0-tz) -g(v_0-tz)-\varepsilon(v_0-tz),
 v_\varepsilon-v_0+tz)\geq 0, \ee
where $t>0$, and $z\in D(L)$ is arbitrary.
Let $\varepsilon\to 0$ in \eqref{e2.5}.
Then, using $v_\varepsilon\rightharpoonup v_0$ and
$Lv_\varepsilon+g(v_\varepsilon)\to 0$,
one gets:
\be \label{e2.6}
 (-L(v_0-tz)-g(v_0-tz), \quad z)\geq 0 \qquad \forall z\in D(L). \ee
Let $t\to 0$ in (2.6). Then $(Lv_0+g(v_0),z)\leq 0\quad\forall z\in 
D(L)$. Since $D(L)$ is dense in $H$, it follows that $Lv_0+g(v_0)=0,$
so $v_0\in N_F$. 

We have proved above that 
$\|v_0\|\leq \|v\|$ for any $v\in N_F$.
Because $N_F$ is a closed and convex set, as we prove below, and
$H$ is a uniformly convex space, there is a 
unique element $v\in N_F$ with minimal norm. Therefore, it 
follows that $v_0=v$, where $v$ is the minimal-norm solution to
(1.2) and $v_0$ is the weak limit of $v_\varepsilon$.

{\it Let us prove the strong convergence $v_\varepsilon\to v$.}

We know that  $v_\varepsilon\rightharpoonup v$.
The inequality $\|v_\varepsilon\|\leq \|v\|$ implies 
$$\|v\|\leq\lim\inf_{\varepsilon\to 0}\|v_\varepsilon\|\leq
\lim\sup_{\varepsilon\to 0}\|v_\varepsilon\|\leq
\|v\|.
$$
Therefore $\lim_{\varepsilon\to 0}\|v_\varepsilon\|=||v||$. Consequently 
one gets:
$$\lim_{\varepsilon\to 0}||v_\varepsilon-v||^2=\lim_{\varepsilon\to 0}
[||v_\varepsilon||^2+||v||^2-2\Re 
(v_\varepsilon,v)]\leq2[||v||^2-(v,v)]=0.$$
Thus,  $v_\varepsilon\to v$, as claimed.

Since $g$ is continuous, it follows that $g(v_\varepsilon)\to g(v)$.

Equation
\be \label{e2.7}
 Lv_\varepsilon +\varepsilon v_\varepsilon + g(v_\varepsilon)=0, \ee
the inequality $\|v_\varepsilon\|\leq \|v\|$, and the relation 
$g(v_\varepsilon)\to g(v)$ imply $Lv_\varepsilon\to \eta:=-g(v)$.

{\it Let us prove that $v_0\in D(L)$}.

Because $v_0=v$, it is sufficient to check
that $v\in D(L)$.  Since $L$ is selfadjoint,
and $Lv_\varepsilon\to \eta$, 
one has:
$$(\eta,\psi)=\lim_{\varepsilon\to 0}(L v_\varepsilon,\psi)
 =\lim_{\varepsilon\to 0}(v_\varepsilon, L\psi)=(v,L\psi) \quad
\forall \psi\in D(L).$$
Thus $v\in D(L)$ and $Lv=\eta$. 

{\it Let us finally check that $N_F$ is closed and convex.}

Assume $F(v_n):=Lv_n+g(v_n)=0$, $v_n\to v$. Then $g(v_n)\to g(v)$
because $g$ is continuous. Thus, $Lv_n\to \eta:=-g(v)$. Since $L$ is 
closed, the relations
$v_n\to v$ and $Lv_n\to \eta$ imply 
$Lv=\eta=-g(v)$. So, $v\in N_F$, and consequently {\it $N_F$ is closed.}

Assume that $0<s<1$, $v,w\in N_F$ and $\psi=sv+(1-s)w$.
Let us show that $\psi\in N_F$. One has $w\in N_F$ if and only if
\be \label{e2.9} (F(z),z-w)\geq 0 \qquad \forall z\in D(L). \ee
Indeed, if $F(w)=0$, then $F(z)=F(z)-F(w)$, and \eqref{e2.9}
holds because $F$ is a monotone operator. 
Conversely, if \eqref{e2.9} holds, then take $z-w=t\eta$, $t=\const>0$,
$\eta\in D(L)$ is arbitrary, and get $F(w+t\eta), \eta)\geq 0$.
Let $t\to 0$, then $(F(w),\eta)\geq 0$ $\forall \eta\in D(L)$.
Thus, since $\overline{D(L)}=H$, one gets $w\in N_F$.
If $\psi=sv+(1-s)w$, and $v,w\in N(F)$,
then $(F(z),z-sv-(1-s)w)=s(F(z),z-v)+(1-s)(F(z),z-w)\geq 0$.
{\it Thus, $N_F$ is convex.}

Theorem 1.2 is proved.      $\Box$

\begin{remark}\label{R:2.2}
\refT{1.2} gives a theoretical framework for a study of nonlinear
ill-posed operator equations
$F(u)=f$, where $F(u)=Lu+g(u)$, where the operator $F'(u)$ is not
boundedly invertible. In this short paper, we do not discuss the case
when the data $f$ are given with some error:
$f_\delta$ is given in place of $f$, $\|f_\delta-f\|\leq \delta$.
In this case one can use DSM for a stable solution of the equation
$F(v)=f$ with noisy data
$f_\delta$, and one uses an algorithm for choosing stopping time
$t_\delta$ such that $\lim_{\delta\to 0} \|u_\delta(t_\delta)-v\|=0$,
where $u_\delta(t)$ is the solution to \eqref{e1.3} with
$\Phi=\Phi_\delta$ chosen suitably (see \cite{R456}).
\end{remark}

\end{document}